\documentclass{birkjour}

\newtheorem*{propa}{Proposition A}

\newcommand{\C}{{\mathbb C}}
\newcommand{\D}{{\mathbb D}}

\newcommand{\T}{{\mathbb T}}

\newcommand{\cN}{{\mathcal N}}

\newcommand{\f}{\frac}

\newcommand{\al}{\alpha}

\newcommand{\ga}{\gamma}

\newcommand{\ze}{\zeta}
\renewcommand{\th}{\theta}

\title[Inner Functions and Inner Factors]
{Inner Functions and Inner Factors\\ 
of Their Derivatives}
\author[K. M. Dyakonov]{Konstantin M. Dyakonov}

\address{%
ICREA and Universitat de Barcelona\\ 
Departament de Matem\`atica Aplicada i An\`alisi\\ 
Gran Via, 585\\ 
E-08007 Barcelona\\ 
Spain}

\email{konstantin.dyakonov@icrea.cat}
\thanks{Supported in part by grant MTM2011-27932-C02-01 from El Ministerio de Ciencia 
e Innovaci\'on (Spain) and grant 2014-SGR-289 from AGAUR (Generalitat de Catalunya).}

\begin{document}
\begin{abstract}
For an inner function $\th$ with $\th'\in\cN$, where $\cN$ is the Nevanlinna class, 
several problems are posed in connection with the canonical (inner-outer) factorization 
of $\th'$. 
\end{abstract}

\maketitle

Let $\D$ stand for the disk $\{z\in\C:|z|<1\}$ and $\T$ for its boundary. Recall that a function $\th\in H^\infty$ 
(i.e., a bounded analytic function $\th$ on $\D$) is called {\it inner} if $\lim_{r\to1^-}|\th(r\ze)|=1$ for 
almost all $\ze\in\T$. Further, a nonvanishing analytic function $F$ on $\D$ is said to be {\it outer} if 
$\log|F|$ agrees with the Poisson integral of an integrable function on $\T$. The Nevanlinna class $\cN$ (resp., 
the Smirnov class $\cN^+$) is formed by the functions representable as $u/v$, where $u,v\in H^\infty$ and 
$v$ is zero-free (resp., outer) on $\D$. Now, for an inner function $\th$, we want to understand what the 
derivative $\th'$ looks like, provided that the latter happens to be in $\cN$ (or $\cN^+$). 

\par We assume that the reader is familiar with basic concepts and facts of function theory on $\D$. These 
include the definitions and standard properties of the Hardy spaces $H^p$, Blaschke products, singular inner 
functions, the canonical factorization in $\cN$, $\cN^+$ and $H^p$, etc. All of this background material can 
be found in \cite[Chapter II]{G}. 

\par The problems raised in this note are in part motivated by the following result; see \cite{DCMFT} or 
\cite{DCR}. 

\begin{propa} Let $\th$ be a nonconstant inner function with $\th'\in\mathcal N$. Then $\th'$ is outer if and 
only if $\th$ is a M\"obius transformation (i.e., a conformal automorphism of $\D$). 
\end{propa}

\par This, in turn, is a simple consequence of the \lq\lq reverse Schwarz--Pick type inequality" 
\begin{equation}\label{eqn:revsp}
\f{1-|\th(z)|^2}{1-|z|^2}\le|\mathcal O(z)|,\qquad z\in\D;
\end{equation}
here $\mathcal O=\mathcal O_{|\th'|}$ is the outer function with modulus $|\th'|$ on $\T$, while $\th$ satisfies 
the hypotheses of Proposition A. For \eqref{eqn:revsp}, the reader is referred to \cite[Section 2]{DSpb} or 
to \cite{DCMFT}, where a more general version is given. 

\par Now, write $\mathfrak I$ for the set of nonconstant inner functions $\th$ with $\th'\in\cN$. For 
$\th\in\mathfrak I$, we actually have the (seemingly) stronger property that $\th'\in\cN^+$, a fact established 
by Ahern and Clark in \cite{AC}. In particular, the map 
\begin{equation}\label{eqn:map}
\th\mapsto\text{\rm inn}(\th'),\qquad\th\in\mathfrak I,
\end{equation}
where $\text{\rm inn}(\cdot)$ stands for \lq\lq the inner factor of", is well defined. Proposition A tells us 
that $J:=\text{\rm inn}(\th')$ is nonconstant for every non-M\"obius $\th\in\mathfrak I$. In addition, the 
{\it boundary spectra} $\sigma(\th)$ and $\sigma(J)$ of the two inner functions must then agree: 
\begin{equation}\label{eqn:sieqsi}
\sigma(\th)=\sigma(J).
\end{equation}
(By definition, the boundary spectrum $\sigma(I)$ of an inner function $I$ is the smallest closed set $K\subset\T$ 
such that $I$ is analytic across $\T\setminus K$.) A weaker version of \eqref{eqn:sieqsi} is verified in \cite{DARX}; 
the full version can be proved by combining inequality \eqref{eqn:revsp} with a result from \cite{KRR}. 

\par Some more notation will be needed. Given a set $\mathfrak F\subset\cN^+$, we write 
$$\text{\rm inn}(\mathfrak F)=\{\text{\rm inn}(f):\,f\in\mathfrak F\},$$
and we denote by $\text{\rm div\,inn}(\mathfrak F)$ the collection of all inner functions that arise as divisors 
of those in $\text{\rm inn}(\mathfrak F)$. Thus, an inner function $I$ is in $\text{\rm div\,inn}(\mathfrak F)$ if 
and only if $J/I\in H^\infty$ for some $J\in\text{\rm inn}(\mathfrak F)$. Clearly, 
$$\text{\rm inn}(\mathfrak F)\subset\text{\rm div\,inn}(\mathfrak F).$$
The set $\mathfrak F$ that chiefly interests us is 
$$\mathfrak I'=\{\th':\,\th\in\mathfrak I\},$$
so that $\text{\rm inn}(\mathfrak I')$ is the range of the map \eqref{eqn:map}. Also, let $\mathfrak B$ (resp., 
$\mathfrak S$) stand for the set of Blaschke products (resp., singular inner functions) lying in $\mathfrak I$. 
We may then look at the classes $\mathfrak B'$ and $\mathfrak S'$, defined as the images of $\mathfrak B$ and 
$\mathfrak S$ under differentiation, and ask about the inner factors (and their divisors) associated with them. 
\par The case of $\mathfrak B'$, however, leads to no new problem. Indeed, for a given $\th\in\mathfrak I$, 
Frostman's theorem (see \cite[Chapter II]{G}) provides us with an $\al\in\D$ such that the function 
$$B_\al:=\f{\th-\al}{1-\bar\al\th}$$ 
is a Blaschke product. Differentiating, we see that $B'_\al\in\cN$ (whence $B_\al\in\mathfrak B$) and 
$\text{\rm inn}(B'_\al)=\text{\rm inn}(\th')$. Therefore, the inner factors that occur for functions in 
$\mathfrak B'$ are the same as those for $\mathfrak I'$. The case of $\mathfrak S'$ is more delicate, though. 

\medskip\noindent\textbf{Problem 1.} Characterize the sets $\text{\rm inn}(\mathfrak I')$ and 
$\text{\rm inn}(\mathfrak S')$. 

\medskip We do not know whether the former set actually coincides with the collection of all inner functions. 
Thus, we ask in particular whether every inner function $J$ can be written as $J=\text{\rm inn}(\th')$ for 
some $\th\in\mathfrak I$. If not, we seek to describe the $J$'s that arise in this way. 
\par As to the other set, $\text{\rm inn}(\mathfrak S')$, it turns out to be {\it strictly} smaller 
than $\text{\rm inn}(\mathfrak I')$ and hence nontrivial. To see why, take $I(z)=zS_1(z)$, where 
\begin{equation}\label{eqn:atomic}
S_1(z):=\exp\left(\f{z+1}{z-1}\right),\qquad z\in\D.
\end{equation}
Letting $a=1-\sqrt2$ and $b(z)=(z-a)/(1-az)$, one verifies by a straightforward calculation that 
$$(bS_1)'(z)=\f{-4I(z)}{(1-az)^2(1-z)^2}$$ 
(in doing so, the identity $a^2-2a-1=0$ should be used). Therefore, $I=\text{\rm inn}((bS_1)')$ and 
hence $I\in\text{\rm inn}(\mathfrak I')$. 
\par On the other hand, if we had $I=\text{\rm inn}(S')$ for some $S\in\mathfrak S$, then \eqref{eqn:sieqsi} 
(with $S$ and $I$ in place of $\th$ and $J$, respectively) would tell us that $\sigma(S)=\sigma(I)=\{1\}$, and 
so $S$ would have to coincide with 
$$S_\ga(z):=\exp\left(\ga\f{z+1}{z-1}\right)$$ 
for some $\ga>0$. However, since $\text{\rm inn}(S'_\ga)=S_\ga$, we see that $I$ does not agree with the 
inner factor of $S'_\ga$ for any $\ga>0$. Consequently, $I\notin\text{\rm inn}(\mathfrak S')$. 

\medskip\noindent\textbf{Problem 2.} Characterize the set $\text{\rm div\,inn}(\mathfrak I')$. 

\medskip While this set contains $\text{\rm inn}(\mathfrak I')$, one may well ask whether the inclusion is 
proper and/or whether every inner function is in $\text{\rm div\,inn}(\mathfrak I')$. 
\par We leave these questions open, but we do observe that $\mathfrak I\subset\text{\rm div\,inn}(\mathfrak I')$. 
Indeed, if $\th\in\mathfrak I$, then $\th$ divides the inner factor of $(\th^2)'=2\th\th'$; this last formula 
also shows that $\th^2$ is in $\mathfrak I$. 
\par Another observation is that, in contrast with the {\it proper} inclusion 
$$\text{\rm inn}(\mathfrak S')\subset\text{\rm inn}(\mathfrak I')$$ 
(see above), we now have 
\begin{equation}\label{eqn:contrast}
\text{\rm div\,inn}(\mathfrak S')=\text{\rm div\,inn}(\mathfrak I').
\end{equation}
In particular, this implies that ${\rm div\,inn}(\mathfrak S')$ is strictly larger than ${\rm inn}(\mathfrak S')$. 

\par To check \eqref{eqn:contrast}, suppose $J$ is a divisor of $\text{\rm inn}(\th')$ for some $\th$ 
in $\mathfrak I$. Now fix a function $S\in\mathfrak S$ (e.g., take $S=S_1$, where $S_1$ is the \lq\lq atomic" 
singular inner function given by \eqref{eqn:atomic}) and put $\varphi:=S\circ\th$. Then $\varphi$ is a singular 
inner function, and since 
$$\varphi'(z)=S'(\th(z))\cdot\th'(z),\qquad z\in\D,$$ 
it follows that $\varphi'\in\cN$ (whence $\varphi\in\mathfrak S$) and that $J$ divides $\text{\rm inn}(\varphi')$. 
This proves the nontrivial inclusion 
$$\text{\rm div\,inn}(\mathfrak S')\supset\text{\rm div\,inn}(\mathfrak I')$$ 
between the two sides of \eqref{eqn:contrast} and thereby establishes the equality. 

\par A similar argument allows us to generalize \eqref{eqn:contrast} as follows. Suppose $E$ is a subset of 
$\D$ such that the class 
$$\mathfrak S_E:=\{\th\in\mathfrak I:\,\th(\D)\cap E=\emptyset\}$$ 
is nonempty. (A result from \cite{AC} tells us that $E$ must be countable.) Then 
\begin{equation}\label{eqn:contrastgen}
\text{\rm div\,inn}(\mathfrak S'_E)=\text{\rm div\,inn}(\mathfrak I').
\end{equation}
Since $\mathfrak S=\mathfrak S_{\{0\}}$, we see that \eqref{eqn:contrast} is indeed a consequence 
of \eqref{eqn:contrastgen}. 

\medskip\noindent\textbf{Problem 3.} For $0<p<1$, let $\mathfrak I^p$ denote the set of all nonconstant 
inner functions $\th$ with $\th'\in H^p$, and let $\mathfrak S^p:=\mathfrak I^p\cap\mathfrak S$. Solve the 
two problems above for $\mathfrak I^p$ and $\mathfrak S^p$ in place of $\mathfrak I$ and $\mathfrak S$. 

\medskip\noindent\textbf{Problem 4.} What, if anything, is the Bergman space counterpart of Proposition A? 
In other words, what happens for \lq\lq Bergman-inner" functions? 

\medskip To be more precise, recall that the Bergman space $A^p$, $0<p<\infty$, is defined as the set 
of analytic functions lying in $L^p(\D)$ with respect to area measure. Recall also that a unit-norm function 
$\psi\in A^p$ is said to be {\it $A^p$-inner} if $|\psi|^p$ annihilates every monomial $z^n$ with $n=1,2,\dots$. 
Furthermore, one defines {\it $A^p$-outer} functions in the appropriate way and proves that every $f\in A^p$ 
has a factorization $f=\psi g$, where $\psi$ is $A^p$-inner and $g$ is $A^p$-outer (see \cite[Chapter 9]{DSch} 
for these notions and results). However, the factorization is not unique. Now, the question is whether 
Proposition A extends to the $A^p$ setting, once suitable adjustments are made. 

\par We conclude with a somewhat vaguer problem, which concerns further possible extensions of Proposition A 
beyond the case of an inner function $\th$. 

\medskip\noindent\textbf{Problem 5.} Find reasonably sharp conditions on a function $f\in H^\infty$ (say, 
with $f'$ in $\cN$ or $\cN^+$) that ensure the presence of a nontrivial inner factor for $f'$, provided 
that $f$ itself enjoys a similar property. Furthermore, assuming that $f$ obeys those conditions and 
letting $I$ and $J$ be the inner factors of $f$ and $f'$, respectively, study the interrelationship 
between $\sigma(I)$ and $\sigma(J)$. 

\medskip Results to that effect could be viewed as descendants of the classical Gauss--Lucas theorem on the 
critical points of a polynomial. In fact, some steps in this direction have already been made in \cite{DARX}. 
The functions $f$ considered there are \lq\lq locally inner" in the sense that $\|f\|_\infty=1$ and $|f|=1$ 
on a set of positive measure (on $\T$). However, these hypotheses seem to be too restrictive, and the results 
in \cite{DARX} are far from being completely satisfactory, so there is plenty of room for improvement and 
further development.

\end{document}